# Continuity in Parametric Linear Programming

by

Somdeb Lahiri

ORCID: https://orcid.org/0000-0002-5247-3497

(Formerly) PD Energy University, Gandhinagar (EU-G), India

somdeb.lahiri@gmail.com

November 1, 2024.

This version: December 8, 2024.

## Abstract

In this paper we assemble some results about the upper-semicontinuity and lower-semicontinuity of the feasible correspondence and the solution correspondence of linear programming problems allowing variability of all parameters of such problems. We also prove continuity properties of optimal value functions, once again allowing all parameters to vary. We discuss sensitivity properties of the optimal value function, keeping the coefficient matrix fixed.

**Keywords:** parametric linear programming, feasibility correspondence, solution correspondence, optimal value function, continuity, upper-semicontinuity, lower-semicontinuity, sensitivity analysis

## 1. Introduction:

Our purpose here is to extend the scope of parametric analysis of linear programming (LP) problems that is available in Holder (2010). Holder (2010) is concerned with (a) sensitivity analysis for the right-hand side of the linear constraints that include in its definition a matrix that remains invariant and (b) sensitivity analysis for the objective function. We intend to allow for variations in all three components of LP problems simultaneously, i.e. the objective function, the right-hand side of the linear constraints and the matrix that determines the left-hand side of the linear constraints. Our main engagement here is to consider continuity of optimal value functions and related results for optimal solution correspondences in the context of LP problems.

The history of parametric linear programming is discussed in Holder (2010), to which we wish to add chapter 8 in the (easy to read) textbook by Saul I Gass (Gass (2003)). Our focus here are the continuity properties of the solution correspondence and the optimal value function of LP problems. The earliest known work in this line of research that we are aware of is Bohm (1975). This was followed by Meyer (1979) and Wets (1985). Our framework of analysis is essentially the one in Wets (1985). The existing literature on continuity properties including this one, uses the concepts of continuity that are available in chapter 1 of Debreu (1959).

After presenting the framework of analysis, we discuss sensitivity analysis for the right-hand side of the constraints. The results we obtain are analogous to the results in the first half of Holder (2010). This is also the framework of analysis that is used in Bohm (1975) and we show that the optimal value function is a concave function of the vector determining the right- hand side of linear constraints. This section is followed by one in which we note that upper-semicontinuity of the feasibility correspondence "comes for free" and prove that if the set of feasible alternatives is bounded for a LP problem, then the feasibility correspondence is lower semicontinuous at that problem. From this section onwards, we rely on continuity concepts defined in chapter 1 of Debreu (1959).

In the next section we provide sufficient conditions for the continuity of the optimal value function and the upper-semicontinuity of the solution function. We show that continuity of the optimal value function implies upper semi-continuity of the solution correspondence. The rest of this section is therefore- and quite naturally- devoted to sufficient conditions for continuity of optimal value functions. A sufficient condition for global continuity of the optimal value function is that the domain is "weakly solution-wise bounded". Another sufficient condition, is the kind of domain considered in Bohm (1975). A sufficient condition for the optimal value function to be continuous at a LP problem is that the problem is regular at that LP problem. The difference between the two results is that the first is concerned with the optimal value function being globally continuous and the second is concerned with continuity of the optimal value function at a point. In this section we also show that global upper-semicontinuity along with a property that we refer to as "weakly solution-wise bounded for the primal" implies global continuity of the optimal value function.

In the subsequent section we discuss the lower-semicontinuity of the solution function, and in the final section of this paper, we replicate the sensitivity analysis in the latter half of Holder (2010) with respect to coefficients determining the objective function. We show that the solution correspondence is lower semicontinuous at a LP problem, if it satisfies the "weakly solution-wise bounded for the primal" property, upper semi-continuity and is singleton-valued at the problem. There are two other results about lower semicontinuity of solution correspondence that may be of some interest. The first result says that if the solution correspondence is singleton valued and regular at a LP problem then the solution correspondence is lower semicontinuous at the problem. The second result says that if a LP problem has a bounded set of feasible points and it is "strongly regular", then the solution correspondence is lower semicontinuous at the problem. For the sensitivity analysis with respect to the objective function of the primal problem, we require that the inequalities in the dual corresponding to the non-basic variables for some basic optimal solution for the primal, hold with strict inequality.

Overall, the properties and domain condition we consider are weaker than the "well-posedness" property discussed in Lucchetti, Radrizzani and Villa (2008).

## 2. Framework of analysis:

Given a positive integer m when we use the term **m-dimensional vector**, we will- unless otherwise specified- mean a point in $\mathbb{R}^m$ say x, whose coordinates are arranged in a column i.e. $x = \begin{pmatrix} x_1 \\ \vdots \\ x_m \end{pmatrix}$.

**Notation:** For vectors x and y of the same dimension, we will use x ≥ y to denote that each coordinate of x is "greater than or equal to" the corresponding coordinate of y, x > y to denote x ≥ y but x ≠ y, and x >> y to denote that each coordinate of x is "strictly greater" than the corresponding coordinate of y.

Given positive integers m and n and m×n **matrix**, say A, is an ordered n-tuple of m-dimensional vectors, whose $j^{th}$ entry for j∈{1, ..., n}, is the m-dimensional column vector denoted by $A^j$.

Thus, $A = (A^1, ..., A^n) \in (\mathbb{R}^m)^n$.

For j∈{1, ..., n}, let $A^j = \begin{pmatrix} a_{1j} \\ \vdots \\ a_{mj} \end{pmatrix}$.

We will denote the set of all m×n matrices as $\mathbb{R}^{m \times n}$ instead of $(\mathbb{R}^m)^n$.

Clearly, an m-dimensional vector is a $\mathbb{R}^{m \times 1}$ matrix.

A 1×n matrix is a **row vector**.

Given an m×n matrix A, its **transpose** denoted by $A^T$ is the n×m matrix whose $i^{th}$ column $(A^T)^i = \begin{pmatrix} a_{i1} \\ \vdots \\ a_{im} \end{pmatrix}$.

For positive integers r, s if B is an r×s matrix with all its columns linearly independent, then $B^T B$ is an s×s matrix which is invertible and the s×r matrix $(B^T B)^{-1} B^T$ denoted $B^-$ is said to be the **pseudo-inverse** of B. Clearly, $B^- B = I^{(s)}$ (the identity matrix of size s) and if B is a square matrix, then $B^- = B^{-1}$.

Given any positive integer r, the **Euclidean norm** of $x \in \mathbb{R}^r$, $\|x\| = \sqrt{\sum_{j=1}^{r}(x_j)^2}$, and for any $A \in \mathbb{R}^{m \times n}$, $\|A\| = \sqrt{\sum_{j=1}^{n} \|A^j\|^2}$.

Given positive integers m, n an n-dimensional vector p, an m-dimensional vector b and an m×n matrix A, we will be concerned with the following m×n **linear programming problem** (LP problem) denoted (p, A, b):

Maximize $p^T x$, subject to Ax = b, x ≥ 0, i.e. $x \in \mathbb{R}_+^n$.

The set of all m×n LP problems is $\mathbb{R}^n \times \mathbb{R}^{m \times n} \times \mathbb{R}^m$. A generic element of $\mathbb{R}^n \times \mathbb{R}^{m \times n} \times \mathbb{R}^m$ will be denoted by ξ and an explicit representation of ξ will often be written as (p(ξ), A(ξ), b(ξ)).

**Important note about notation used in what follows:** Given a sequence <ξ(N)|N∈ℕ> in $\mathbb{R}^n \times \mathbb{R}^{m \times n} \times \mathbb{R}^m$, to economize on notational complexity, we will refer to (p(ξ(N)), A(ξ(N)), b(ξ(N))) as (p(N), A(N), b(N)) for all N∈ℕ, whenever there is no scope for confusion.

Given a sequence <ξ(N)|N∈ℕ> in $\mathbb{R}^n \times \mathbb{R}^{m \times n} \times \mathbb{R}^m$ and $\xi = (p, A, b) \in \mathbb{R}^n \times \mathbb{R}^{m \times n} \times \mathbb{R}^m$ we will say that <ξ(N)|N∈ℕ> **converges to** (p, A, b) or ξ and write $\lim_{N \to \infty} \xi(N) = (p, A, b)$ to mean the following:

(i) $\lim_{N\to\infty} ||p(N) - p|| = 0$ and $\lim_{N\to\infty} ||b(N) - b|| = 0$.

(ii) $\lim_{N\to\infty} ||A(N) - A|| = 0$.

Given, $\xi = (p, A, b) \in \mathbb{R}^n \times \mathbb{R}^{m\times n} \times \mathbb{R}^m$, let $X(p, A, b) = X(\xi) = \{x \in \mathbb{R}_+^n | Ax = b\}$ and $X^*(\xi) = \{x \in X(\xi) | \{A^j | x_j > 0\}$ is linearly independent$\}$.

Suppose that for $\xi \in \mathbb{R}^n \times \mathbb{R}^{m\times n} \times \mathbb{R}^m$, $X(\xi)$ is non-empty. Then it is a convex set and $X^*(\xi)$ is a non-empty and finite set (see Lahiri (2020)). If in addition $X(\xi)$ is bounded, then Corollary of theorem 3 in Lahiri (2024), says that $X(\xi)$ is the convex hull of $X^*(\xi)$ and points in $X^*(\xi)$ are "extreme points" of $X(\xi)$.

Given, $\xi = (p, A, b) \in \mathbb{R}^n \times \mathbb{R}^{m\times n} \times \mathbb{R}^m$, let $S(\xi) = S(p, A, b) = \{x \in X(p, A, b) | p^T x \geq p^T x'$ for all $x' \in X(p, A, b)\}$.

$S(\xi)$ is said to be **the set of optimal solutions** at $\xi$.

It is well known that if $S(\xi) \neq \phi$ then $S^*(\xi) = S(\xi) \cap X^*(\xi) \neq \phi$ (see Lahiri (2020)).

If $X(\xi)$ is non-empty and bounded, then $S(\xi) \neq \phi$ (see Lahiri (2020)), and from theorem 3 in Lahiri (2024), we know that: (i) $S(\xi)$ is the convex hull of $S^*(\xi)$, (ii) every point in $S^*(\xi)$ is an extreme point of $S(\xi)$.

Given $\xi \in \mathbb{R}^n \times \mathbb{R}^{m\times n} \times \mathbb{R}^m$, its dual denoted $D(\xi)$ is the following minimization problem.

Minimize $y^T b(\xi)$, subject to $y^T A(\xi) \geq p(\xi)^T$, y unconstrained in sign.

The following results are well known and available with all the required proofs in Lahiri (2020).

(i) $x \in \mathbb{R}^n$ solves $\xi$ <u>if and only if</u> $x \in X(\xi)$ and there exists $y \in \mathbb{R}^m$ satisfying $y^T A(\xi) \geq p(\xi)^T$, $(y^T A(\xi) - p(\xi)^T)x = 0$.

(ii) Suppose $x \in S^*(\xi)$ and $B(\xi) = \{A^j(\xi) | x_j > 0\}$ where the columns in $B(\xi)$ are linearly independent. Suppose that (by re-arranging the co-ordinates of $A(\xi)$ if necessary) $A(\xi) = [B(\xi)|E(\xi)]$ so that $x = \begin{pmatrix} x_{B(\xi)} \\ x_{E(\xi)} \end{pmatrix} = \begin{pmatrix} x_{B(\xi)} \\ 0 \end{pmatrix}$ and $p = \begin{pmatrix} p_{B(\xi)} \\ p_{E(\xi)} \end{pmatrix}$.

Then $x_{B(\xi)} = (B(\xi)^T B(\xi))^{-1} B(\xi)^T b(\xi) = B(\xi)^- b(\xi)$ and y such that $y^T = (p_{B(\xi)})^T B(\xi)^-$ solves $D(\xi)$. Further $p(\xi)^T x = (p_{B(\xi)})^T x_{B(\xi)} = y^T b(\xi)$.

(i) above is known as the **Karush-Kuhn-Tucker (KKT) conditions** for $\xi$.

**Note:** While $B(\xi)$ and $E(\xi)$ are arrays of columns in $A(\xi)$, since the columns in $B(\xi)$ are linearly independent, they must be distinct. Thus, $B(\xi) = <A^j(\xi) | x_j > 0>$ can also be written as $\{A^j(\xi) | x_j > 0\}$. However, the same is <u>not applicable</u> for $E(\xi) = <A^j(\xi) | x_j = 0>$, since it is not necessary that all columns of $E(\xi)$ are distinct.

It is easy to see that if $<\xi(N) | N \in \mathbb{N}>$ **converges to** $(p, A, b)$, $A = [B|E]$ with the columns of B being linearly independent and $A(\xi(N)) = [B(\xi(N)) | E(\xi(N))]$ for all $N \in \mathbb{N}$, with the columns of B(ξ(N)) being linearly independent for all N ≥ N₁, for some N₁∈ ℕ, then the sequence of matrices <C(N)|N∈ℕ> with C(N) = B⁻(ξ(N)) for all N ≥ N₁ converges to B⁻.

Let $\mathcal{X}$ be a non-empty subset of $\mathbb{R}^n \times \mathbb{R}^{m \times n} \times \mathbb{R}^m$, such that every for every ξ∈$\mathcal{X}$, S(ξ) is non-empty. Thus, for all ξ∈$\mathcal{X}$, $S^*(\xi) \neq \phi$.

$\mathcal{X}$ is said to be **the domain**.

Given (p, A, b)∈$\mathcal{X}$ and δ > 0, a δ-**neighborhood** of (p, A, b) is the set {ξ∈ $\mathcal{X}$| ||p(ξ) − p|| < δ, ||A(ξ) − A|| < δ and ||b(ξ) − b|| < δ}.

The function V: $\mathcal{X} \to \mathbb{R}$ such that for all ξ∈$\mathcal{X}$, V(ξ) = p(ξ)ᵀx, x∈S(ξ), is said to be the **optimal value function** on $\mathcal{X}$.

Thus, for all ξ∈$\mathcal{X}$, S(ξ) = {x∈$\mathbb{R}^n_+$| A(ξ)x = b(ξ), p(ξ)ᵀx ≥ V(ξ)} = {x∈$\mathbb{R}^n_+$| A(ξ)x = b(ξ), p(ξ)ᵀx = V(ξ)}.

The following concepts are available in Debreu (1959).

For r, s and a non-empty subsets X of $\mathbb{R}^r$ and Y of $\mathbb{R}^s$, a function F: X→ $2^Y$(i.e., the power set of Y) is said to be a correspondence from X to Y and is denoted by F:X→→Y.

Given x∈X, F is said to be **upper semicontinuous at** x, if for all sequences <x(N)|N∈ℕ> in X converging to x∈X and all sequences <y(N)|N∈ℕ> in Y converging y∈Y: [y(N)∈F(x(N)) for all N∈ℕ] implies [y∈F(x)].

F is said to be **upper semicontinuous (on X)** if it is upper semicontinuous at all x∈X.

Given x∈X, F is said to be **lower semicontinuous at** x, if for all sequences <x(N)|N∈ℕ> in X converging to x∈X and y∈F(x), there exists a sequence <y(N)|N∈ℕ> in Y converging to y that satisfies y(N)∈F(x(N)) for all N∈ℕ.

F is said to be **lower semicontinuous (on X)** if it is lower semicontinuous at all x∈X.

Let X: $\mathcal{X} \to \to \mathbb{R}^n_+$ be the correspondence that assigns to each (p, A, b)∈$\mathcal{X}$ the set X(p, A, b) = {x∈$\mathbb{R}^n_+$| Ax = b}. We will refer to X as the **feasible correspondence** on $\mathcal{X}$.

Let S:$\mathcal{X} \to \to \mathbb{R}^n_+$ be the correspondence that assigns to each ξ∈$\mathcal{X}$ the set S(ξ) defined above. We will refer to S as the **solution correspondence** on $\mathcal{X}$.

### 3. Feasibly convex domain and the concavity of the optimal value function:

$\mathcal{X}$ is said to be **feasibly convex** if (i) there exists A∈$\mathbb{R}^{m \times n}$ and p∈$\mathbb{R}^n$, such that for all ξ∈$\mathcal{X}$, A(ξ) = A and p(ξ) = p; (ii) (p, A, $b^{(1)}$), (p, A, $b^{(2)}$) ∈$\mathcal{X}$ and t∈[0, 1] implies (p, A, $tb^{(1)}$ + (1-t)$b^{(2)}$)∈$\mathcal{X}$.

That there is an m×n matrix A and an n-vector p, such that for all ξ∈$\mathcal{X}$, A(ξ) = A and p(ξ) = p, is explicitly assumed for the entire analysis in Bohm (1975).

Holder (2010) discusses the sensitivity analysis under the domain assumption in Bohm (1975).

Given, (p, A, b)∈ $\mathcal{X}$, let Δb be a non-zero m-dimensional column vector.

Let $x = \begin{pmatrix} x_B \\ x_E \end{pmatrix} = \begin{pmatrix} x_B \\ 0 \end{pmatrix} \in S^*(p, A, b)$, where $A = (B|E)$ with the columns of A in the submatrix B being linearly independent and let $p = \begin{pmatrix} p_B \\ p_E \end{pmatrix}$. Thus, $x_B = B^-b \gg 0$.

Thus, $y^T = (p_B)^T B^-$ satisfies $y^T B = (p_B)^T$ and $y^T E \geq (p_E)^T$.

Consider the LP problem $(p, A, b + \theta \Delta b)$.

Since, $y^T = (p_B)^T B^-$ satisfies $y^T B = (p_B)^T$ and $y^T E \geq (p_E)^T$, $\begin{pmatrix} B^-(b + \theta \Delta b) \\ 0 \end{pmatrix}$ along with y satisfy the KKT conditions for $(p, A, b + \theta \Delta b)$, so long as so long as $B^-(b + \theta \Delta b) \geq 0$, i.e., $[B^-]_i(b + \theta \Delta b) \geq 0$ for all i corresponding to the $i^{th}$ coordinate of $x_B$. Thus, $\begin{pmatrix} B^-(b + \theta \Delta b) \\ 0 \end{pmatrix}$ is an optimal solution for $(p, A, b + \theta \Delta b)$ so long as $B^-(b + \theta \Delta b) \geq 0$.

Thus, for $\theta$ satisfying $\theta \geq \max \{- \frac{[B^-]_i b}{[B^-]_i \Delta b} | [B^-]_i \Delta b > 0\}$, and $\theta \leq \min \{- \frac{[B^-]_i b}{[B^-]_i \Delta b} | [B^-]_i \Delta b < 0\}$ it must be the case that $\begin{pmatrix} B^-(b + \theta \Delta b) \\ 0 \end{pmatrix}$ is an optimal solution for $(p, A, b + \theta \Delta b)$.

Thus, a sufficient condition for $\begin{pmatrix} B^-(b + \theta \Delta b) \\ 0 \end{pmatrix}$ to be an optimal solution for $(p, A, b + \theta \Delta b)$ is that $\theta$ satisfies the following inequality: $\max \{- \frac{[B^-]_i b}{[B^-]_i \Delta b} | [B^-]_i \Delta b > 0\} \leq \theta \leq \min \{- \frac{[B^-]_i b}{[B^-]_i \Delta b} | [B^-]_i \Delta b > 0\}$.

Since $\max \{- \frac{[B^-]_i b}{[B^-]_i \Delta b} | [B^-]_i \Delta b > 0\} < 0 < \min \{- \frac{[B^-]_i b}{[B^-]_i \Delta b} | [B^-]_i \Delta b > 0\}$, there is a non-degenerate closed interval around zero such that a sufficient condition for $\begin{pmatrix} B^-(b + \theta \Delta b) \\ 0 \end{pmatrix}$ to be an optimal solution for $(p, A, b + \theta \Delta b)$ is that $\theta$ lies in this interval.

For all values of $\theta$ in this interval $V(p, A, b + \theta \Delta b) = p_B^T B^-(b + \theta \Delta b) = V(p, A, b) + \theta p_B^T B^- \Delta b$.

Hence in this interval V is a differentiable function of $\theta$ with its constant derivative (slope) being $p_B^T B^- \Delta b$.

It is important to note that this linearity condition of the optimal value function does not depend on whether for $\theta \neq 0$, $(p, A, b + \theta \Delta b) \in \mathcal{X}$ or not.

**Proposition 1:** If $\mathcal{X}$ is feasibly convex, then V is concave.

**Proof:** Let $(p, A, b^{(1)}), (p, A, b^{(2)}) \in \mathcal{X}$ and $t \in (0, 1)$. Since $\mathcal{X}$ is convex, $(p, A, tb^{(1)} + (1-t)b^{(2)}) \in \mathcal{X}$.

Let $x^{(k)} \in S(p, A, b^{(k)})$ for $k \in \{1, 2\}$. Thus, $x^{(k)} \geq 0$ and $Ax^{(k)} = b^{(k)}$.

Further $V(p, A, b^{(k)}) = p^T x^{(k)}$.

Clearly, $tx^{(1)} + (1-t)x^{(2)} \geq 0$ and $A(tx^{(1)} + (1-t)x^{(2)}) = tb^{(1)} + (1-t)b^{(2)}$.

Thus, $tx^{(1)} + (1-t)x^{(2)} \in X(p, A, tb^{(1)} + (1-t)b^{(2)})$.

Thus, $V(p, A, tb^{(1)} + (1-t)b^{(2)}) \geq p^T(tx^{(1)} + (1-t)x^{(2)}) = tp^Tx^{(1)} + (1-t)px^{(2)} = tV(p, A, b^{(1)}) + (1-t)V(p, A, b^{(2)})$. Q.E.D.

## 4. Upper semicontinuity and Lower-semicontinuity of the feasibility correspondence:

**Proposition 2:** (a) The feasible correspondence X on $\mathcal{X}$ is upper semicontinuous. (b) If at $\xi \in \mathcal{X}$, $X(\xi)$ is bounded then the correspondence X is lower semicontinuous <u>at</u> $\xi$.

**Proof:** The proof of part (a) is quite straightforward, hence let us prove part (b).

Suppose that at $\xi = (p, A, b)$, $X(\xi)$ is bounded at $\xi \in \mathcal{X}$. Then, since $X(\xi)$ is non-empty as well, Corollary of theorem 3 in Lahiri (2024) says that $X(\xi)$ is the convex hull of $X^*(\xi) = \{x \in X(\xi) | \{A^j | x_j > 0\}$ is linearly independent$\}$ and points in $X^*(\xi)$ are "extreme points" of $X(\xi)$.

Thus, if we can show that for each $x \in X^*(\xi)$- the latter being a finite set- and any sequence $<\xi(N)|N \in \mathbb{N}>$, there exists a sequence $<x(N)|N \in \mathbb{N}>$ converging to x and satisfies $x(N) \in X(\xi(N))$ for all $N \in \mathbb{N}$, then we are done.

Hence, let $x \in X^*(\xi)$ and without loss of generality suppose that $\{A^j | x_j > 0\}$ precede all other columns of A, so that $A = (B|E)$ with $B = \{A^j | x_j > 0\}$ and the columns of B are linearly independent. Let $x = \begin{pmatrix} x_B \\ x_E \end{pmatrix} = \begin{pmatrix} x_B \\ 0 \end{pmatrix} = \begin{pmatrix} B^-b \\ 0 \end{pmatrix}$ and $B^-b \gg 0$.

Since $B^-b \gg 0$ and since $<(A(N), b(N))|N \in \mathbb{N}>$ converges to $(A, b)$ there exists $N_1 \in \mathbb{N}$, such that for all $N \geq N_1$, the columns in $B(N) = \{A^j(N) | x_j > 0\}$ are linearly independent and $B(N)^-b(N) \gg 0$. Thus, $x(N) = \begin{pmatrix} x_{B(N)} \\ x_{E(N)} \end{pmatrix} = \begin{pmatrix} x_{B(N)} \\ 0 \end{pmatrix} = \begin{pmatrix} B(N)^-b(N) \\ 0 \end{pmatrix} \in X(\xi(N))$ for all $N \geq N_1$.

Let $x(N) \in X(\xi(N))$ for all $N < N_1$ and $x(N) = \begin{pmatrix} B(N)^-b(N) \\ 0 \end{pmatrix}$ for all $N \geq N_1$.

Then, the sequence $<x(N)|N \in \mathbb{N}>$ converges to x. Q.E.D.

## 5. Upper-semicontinuity of solution correspondence and continuity of optimal value function:

The proof of the following result parallels the one in Debreu (1959) about the upper-semicontinuity of the demand correspondence.

**Proposition 3:** If at $\xi \in \mathcal{X}$, $X(\xi)$ is bounded then the correspondence S is upper semicontinuous <u>at</u> $\xi$.

**Proof:** Let $\xi = (p, A, b) \in \mathcal{X}$ and suppose $X(\xi)$ is bounded. Suppose $<\xi(N)|N \in \mathbb{N}>$ is a sequence in $\mathcal{X}$ converging to $\xi$.

Let $<x(N)|N \in \mathbb{N}>$ be a sequence in $\mathbb{R}^n$ such that $x(N) \in S(\xi(N))$ for all $N \in \mathbb{N}$ and suppose $<x(N)|N \in \mathbb{N}>$ converges to x.

Towards a contradiction suppose $x \notin S(\xi)$. Thus, there exists $x^* \in X(\xi)$ such that $p^T x^* > p^T x$.

Since X(ξ) is bounded, by part (b) of proposition 2, the feasibility correspondence is lower semicontinuous at ξ.

Thus, for each N∈ℕ, there exists $x^*(N) \in X(\xi(N))$ such that the sequence $\langle x^*(N) | N \in \mathbb{N} \rangle$ converges to $x^*$.

Since $p^T x^* > p^T x$, $\langle x(N) | N \in \mathbb{N} \rangle$ converges to x and $\langle x^*(N) | N \in \mathbb{N} \rangle$ converges to $x^*$ implies that there exists $N_1 \in \mathbb{N}$ such that for all $N \geq N_1$, $p^T x^*(N) > p^T x(N)$.

However, $x^*(N) \in X(\xi(N))$ for all $N \in \mathbb{N}$, and $p^T x^*(N) > p^T x(N)$, contradicts $x(N) \in S(\xi(N))$ for all $N \in \mathbb{N}$.

Thus, it must be the case that $x \in S(\xi)$, i.e. S is upper semicontinuous at ξ. Q.E.D.

Our next result is an important reason for whatever follows it in this section.

**Proposition 4:** If V is continuous on $\mathcal{X}$, then the solution correspondence S on $\mathcal{X}$ is upper semicontinuous.

**Proof:** Let ξ = (p, A, b) ∈ $\mathcal{X}$ and suppose $\langle \xi(N) | N \in \mathbb{N} \rangle$ is a sequence in $\mathcal{X}$ converging to ξ. Let $\langle x(N) | N \in \mathbb{N} \rangle$ be a sequence in $\mathbb{R}^n$ such that $x(N) \in S(\xi(N))$ for all $N \in \mathbb{N}$ and suppose $\langle x(N) | N \in \mathbb{N} \rangle$ converges to x. Thus, $V(\xi(N)) = p(N)^T x(N)$ for all $N \in \mathbb{N}$.

Since $\langle \xi(N) | N \in \mathbb{N} \rangle$ converges to ξ by the continuity of V, $\langle V(\xi(N)) | N \in \mathbb{N} \rangle$ converges to V(ξ).

However, $\langle \xi(N) | N \in \mathbb{N} \rangle$ converges to ξ and $\langle x(N) | N \in \mathbb{N} \rangle$ converges to x implies $\langle p(N)^T x(N) | N \in \mathbb{N} \rangle$ converges to $p^T x$.

Thus, $V(\xi) = p^T x$.

$\langle \xi(N) | N \in \mathbb{N} \rangle$ converges to ξ and $\langle x(N) | N \in \mathbb{N} \rangle$ converges to x implies $\langle A(N)x(N) | N \in \mathbb{N} \rangle$ converges to Ax and $\langle b(N) | N \in \mathbb{N} \rangle$ converges to b.

Since A(N)x(N) = b(N) for all N∈ℕ, taking limits on both sides we get Ax = b.

Further, $\langle x(N) | N \in \mathbb{N} \rangle$ converges to x and $x(N) \geq 0$ for all $N \in \mathbb{N}$ implies $x \geq 0$.

Thus, Ax = b, x ≥ 0 and $V(\xi) = p^T x$, whence $x \in S(\xi)$.

Thus, S is upper semicontinuous. Q.E.D.

Our framework of analysis being very general we need additional assumptions to obtain our desired results.

Recall that for $\xi = (p, A, b) \in \mathbb{R}^n \times \mathbb{R}^{m \times n} \times \mathbb{R}^m$ and any pair $(x,y) \in \mathbb{R}^n \times \mathbb{R}^m$: [x solves ξ and y solves D(ξ)] <u>if and only if</u> [Ax = b, x ≥ 0, $y^T A \geq p^T$ and $(y^T A - p^T)x = 0$].

$\mathcal{X}$ is said to be **weakly solution-wise bounded** if there is a non-empty bounded subset Γ of $\mathbb{R}^n \times \mathbb{R}^m$ such that for all $\xi \in \mathcal{X}$: $\Gamma \cap \{(x,y) \in \mathbb{R}^n \times \mathbb{R}^m | x \in X(\xi), y^T A(\xi) \geq p(\xi)^T$ and $(y^T A(\xi) - p(\xi)^T)x = 0\} \neq \phi$.

**Proposition 5:** Suppose $\mathcal{X}$ is weakly solution-wise bounded. Then V is continuous on $\mathcal{X}$.

**Proof:** Let $\xi = (p, A, b) \in \mathcal{X}$ and suppose $\langle \xi(N) | N \in \mathbb{N} \rangle$ is a sequence in $\mathcal{X}$ converging to $\xi$. Towards a contradiction suppose there exists $\varepsilon > 0$ such that $|V(\xi(N)) - V(\xi)| \geq \varepsilon$ infinitely often. Without loss of rigor or compromising the validity of our proof, suppose $|V(\xi(N)) - V(\xi)| \geq \varepsilon$ for all $N \in \mathbb{N}$. Such an assumption helps to economize on notation.

By our assumption that $\mathcal{X}$ is weakly solution-wise bounded, for each $N \in \mathbb{N}$, there exists a bounded subset $\Gamma$ of $\mathbb{R}^n \times \mathbb{R}^m$ and $(x(N), y(N)) \in \Gamma$ such that $A(N)x(N) = b(N)$, $x(N) \geq 0$, $y(N)^T A(N) \geq p(N)^T$, $(y(N)^T A(N) - p(N)^T)x(N) = 0$ for all $N \in \mathbb{N}$.

The sequence $\langle (x(N), y(N)) | N \in \mathbb{N} \rangle$ in $\Gamma$ is bounded since $\Gamma$ is a non-empty bounded subset of $\mathbb{R}^n \times \mathbb{R}^m$.

By the Bolzano-Weierstrass's theorem, (a proof of which is available in Topic 5 of Lahiri (2020)) we know that $\langle (x(N), y(N)) | N \in \mathbb{N} \rangle$ has a subsequence $\langle (x(k(N)), y(k(N))) | N \in \mathbb{N} \rangle$ converging to $(x, y) \in \mathbb{R}^n \times \mathbb{R}^m$.

Further, $|V(\xi(N)) - V(\xi)| \geq \varepsilon$ for all $N \in \mathbb{N}$ implies $|V(\xi(k(N))) - V(\xi)| \geq \varepsilon$ for all $N \in \mathbb{N}$.

Since, (i) $\langle (p(N), A(N), b(N)) | N \in \mathbb{N} \rangle$ converges to $(p, A, b)$, (ii) $\langle (x(k(N)), y(k(N))) | N \in \mathbb{N} \rangle$ converges to $(x, y)$, and (iii) $[A(k(N))x(k(N)) = b(k(N)), x(k(N)) \geq 0, y(k(N))^T A(k(N)) \geq p(k(N))^T, (y(k(N))^T A(k(N)) - p(k(N))^T)x(k(N)) = 0$ for all $N \in \mathbb{N}]$, we get $Ax = \lim_{N \to \infty} A(k(N))x(k(N)) = \lim_{N \to \infty} b(k(N)) = b$, $x = \lim_{N \to \infty} x(k(N)) \geq 0$, $y^T A = \lim_{N \to \infty} y(k(N))^T A(k(N)) \geq \lim_{N \to \infty} p(k(N))^T = p^T$, $0 = \lim_{N \to \infty} (y(k(N))^T A(k(N)) - p(k(N))^T)x(k(N)) = (y^T A - p^T)x$.

Since $Ax = b$, $x \geq 0$, $y^T A \geq p^T$ and $(y^T A - p^T)x = 0$, $x \in S(\xi)$ and $y$ solves $D(\xi)$ and hence $V(\xi) = p^T \xi$.

Further, $V(\xi(N)) = p(N)^T x(N)$ for all $N \in \mathbb{N}$.

Thus, $V(\xi) = p^T x = \lim_{N \to \infty} p(k(N))^T x(k(N)) = \lim_{N \to \infty} V(\xi(k(N))$, contradicting $V(\xi(k(N))) - V(\xi)| \geq \varepsilon$ for all $N \in \mathbb{N}$.

Thus, it must be the case that $\langle V(\xi(N)) | N \in \mathbb{N} \rangle$ converges to $V(\xi)$.

This being true for all $\xi \in \mathcal{X}$, we may conclude that $V$ is continuous on $\mathcal{X}$. Q.E.D.

**Note:** Our assumption that $\mathcal{X}$ is said to be "weakly solution-wise bounded" is required for the validity of Proposition 5 as the following example reveals.

**Example 1:** For $N \in \mathbb{N}$, let $\xi(N)$ be the LP problem:

Maximize $\frac{1}{N} x_1$, subject to $\frac{1}{N} x_1 + x_2 = 1$, $x_1, x_2 \geq 0$.

$S(\xi(N)) = \{(N, 0)\}$ and $V(\xi(N)) = 1$ for all $N \in \mathbb{N}$.

The sequence $\langle \xi(N) | N \in \mathbb{N} \rangle$ converges to $(\binom{0}{0}, (0|1), (1))$, $S(\binom{0}{0}, (0|1), (1)) = \{(x_1, 1) | x_1 \geq 0\}$ and $V(\binom{0}{0}, (0|1), (1)) = 0$. Thus, $V$ is not continuous at $(\binom{0}{0}, (0|1), (1))$.

An alternative away of approaching the continuity of optimal value function could be this.

$(p, A, b) \in \mathbb{R}^n \times \mathbb{R}^{m \times n} \times \mathbb{R}^m$ is said to be **regular** if there exists $x \in S^*(p, A, b)$ such that for $y^T = (p_B)^T B^-$ where B is the submatrix of A whose columns are $\{A^j | x_j > 0\}$ and $p_B$ is the sub-vector of p that correspond to the columns of B, it is the case that $y^T A^j - p_j > 0$ for all $j \in \{1, \ldots, n\}$ satisfying $x_j = 0$.

**Proposition 6:** Suppose $\xi \in \mathcal{X}$ is regular. Then, V is continuous at $\xi$.

**Proof:** As in the proof of Proposition 5, suppose $<\xi(N)|N \in \mathbb{N}>$ is a sequence in $\mathcal{X}$ converging to $\xi = (p, A, b)$.

Since $\xi$ is regular, there exists $x \in S^*(\xi)$ and for $y^T = (p_B)^T B^-$ where $B = \{A^j | x_j > 0\}$ it is the case that $y^T A^j - p_j > 0$ for all $j \in \{1, \ldots, n\}$ satisfying $x_j = 0$.

Clearly, $y^T A^j = p_j$ for all $j \in \{1, \ldots, n\}$ satisfying $x_j = 0$, so that $y^T A \geq p^T$ and $(y^T A - p^T)x = 0$

By re-arranging the co-ordinates of A if necessary, let $A = [B|E]$.

Thus, $x = \begin{pmatrix} x_B \\ x_E \end{pmatrix} = \begin{pmatrix} x_B \\ 0 \end{pmatrix}$ satisfies $x_B = B^- b \gg 0$.

Let $p = \begin{pmatrix} p_B \\ p_E \end{pmatrix}$. Thus, $V(\xi) = p^T x = (p_B)^T x_B = (p_B)^T B^- b$.

$y^T A^j - p_j = 0$ for all $j \in \{1, \ldots, n\}$ satisfying $x_j > 0$ implies $y^T B = (p_B)^T$.

Since, $y^T A \geq p^T$, $(y^T A - p^T)x = 0$ and $y^T B = (p_B)^T$, from Topic 2 of Lahiri (2020) we know that $y^T = (p_B)^T B^-$ so that $y^T b = (p_B)^T B^- b = (p_B)^T x_B = p^T x$.

Further, $y^T A^j - p_j > 0$ for all $j \in \{1, \ldots, n\}$ satisfying $x_j = 0$ implies $(p_B)^T B^- A^j > p_j$ for all $j \in \{1, \ldots, n\}$ satisfying $x_j = 0$.

Since, $<\xi(N)|N \in \mathbb{N}>$ is a sequence in $\mathcal{X}$ converging to $\xi$ there exists $N_1 \in \mathbb{N}$ such that for all $N \geq N_1$, the columns of B(N) are linearly independent, $x(N) = \begin{pmatrix} x_{B(N)} \\ 0 \end{pmatrix} = \begin{pmatrix} B(N)^- b(N) \\ 0 \end{pmatrix}$ along with $y(N)^T = (p_{B(N)})^T B(N)^-$ satisfies $A(N)x(N) = b$, $x_j(N) > 0$ if and only if the $j^{th}$ column of A is included in B(N), $y(N)^T B(N) = (p_{B(N)})^T$ and $y(N)^T A^j(N) > p_j(N)$ for all j such that $x_j(N) = 0$.

Thus, $V(\xi(N)) = p(N)^T x(N) = (p_B(N))^T x_B(N)$ for all $N \geq N_1$.

Let $<x(N) N \in \mathbb{N}>$ be such that $x(N) \in S(\xi(N))$ for all $N < N_1$ and $x(N) = \begin{pmatrix} B(N)^- b(N) \\ 0 \end{pmatrix}$ for all $N \geq N_1$.

Thus, $V(\xi(N)) = p(N)^T x(N)$ for all $N \in \mathbb{N}$

Since $<(p(N), A(N), b(N)|N \in \mathbb{N}>$ converges to $(p, A, b)$, it must be the case that $<(p_B(N), b(N), b(N)|N \in \mathbb{N}>$ converges to $(p_B, B, b)$.

Thus, $<x(N)|N \in \mathbb{N}>$ converges to $x = \begin{pmatrix} x_B \\ 0 \end{pmatrix} = \begin{pmatrix} B^- b \\ 0 \end{pmatrix}$ and hence $<p(N)^T x(N)|N \in \mathbb{N}>$ converges to $p_B B^- b$.

Thus, $\langle V(\xi(N)) | N \in \mathbb{N} \rangle$ converges to $V(\xi)$. Hence V is continuous at $\xi$. Q.E.D.

As in Bohm (1975), if we had required that there is an m×n matrix A and an n-vector p, such that for all $\xi \in \mathcal{X}$, $A(\xi) = A$ and $p(\xi) = p$, then the assumptions in propositions 5 and 6 can be dispensed with.

The following is the equivalent of theorem 1 in Bohm (1975), whose proof is being provided for the following reasons: (i) The formulation of LP problems in Bohm (1975) is different from our formulation of LP problems. (ii) Bohm (1975) does not prove the result and refers to Wetts (1966) for a proof. (iii) In comparison to our work here, the analysis in Bohm (1975) is concerned with a very special case.

**Proposition 7:** Suppose that there is an m×n matrix A and an n-vector p, such that for all $\xi \in \mathcal{X}$, $A(\xi) = A$ and $p(\xi) = p$. Then V is continuous on $\mathcal{X}$.

**Proof:** Let $\xi \in \mathcal{X}$ and suppose $\langle \xi(N) | N \in \mathbb{N} \rangle$ is a sequence in $\mathcal{X}$ converging to $\xi$.

To economize on notational complexity we will refer to $b(\xi)$ as b and to $b(\xi(N))$ as b(N) for all $N \in \mathbb{N}$.

There exists $x \in S^*(\xi)$, so that (by re-arranging the co-ordinates of A if necessary) A = [B|E], where the columns of B are linearly independent, $x = \binom{x_B}{x_E} = \binom{x_B}{0}$, $p = \binom{p_B}{p_E}$ satisfying (i) $x_B$ = B⁻b with all its coordinates being strictly positive and (ii) $y = (p_B)^T B^-$ satisfies $y^T B = (p_B)^T$ and $y^T A^j \geq p_j$ for all j such that $x_j = 0$. Thus, $y^T A \geq p^T$ and $(y^T A - p^T)x = 0$.

Thus, $V(\xi) = p^T x = (p_B)^T x_B = (p_B)^T B^- b$.

Since, $\langle \xi(N) | N \in \mathbb{N} \rangle$ is a sequence in $\mathcal{X}$ converging to $\xi$ there exists $N_1 \in \mathbb{N}$ such that for all $N \geq N_1$, such that $x(N) = \binom{B(N)^- b(N)}{0}$ with $x_j(N) > 0$ if and only if the j$^{th}$ column of A is included in B, $y^T A \geq p^T$ and $(y^T A - p^T)x = 0$.

Thus, $V(\xi(N)) = p^T x(N) = (p_B)^T B(N)^- b(N)$ for all $N \geq N_1$.

Let $x(N) \in S(\xi(N))$ for all $N < N_1$ and $x(N) = \binom{B(N)^- b(N)}{0}$ for all $N \geq N_1$.

Since $\langle b(N) | N \in \mathbb{N} \rangle$ converges to b, $\langle x(N) | N \in \mathbb{N} \rangle$ converges to B⁻b and thus $\langle V(\xi(N)) | N \in \mathbb{N} \rangle$ converges to $V(\xi)$.

Thus, V is continuous at $\xi$. $\xi \in \mathcal{X}$ being arbitrary, we get the desired result. Q.E.D.

**Note:** Theorem 2 of Wets (1985) says that if the correspondence from $\mathcal{X}$ to $\mathbb{R}^n \times \mathbb{R}^m$ defined by $\xi \mapsto X(\xi) \times \{y \in \mathbb{R}^m | y^T A(\xi) \geq p(\xi)^T\}$ is both upper and lower continuous at an LP problem, then (under assumptions implicit in the analysis) V is continuous at the same LP problem..

We will now show that along with a "weaker version" of the weakly solution-wise bounded property, the upper-semicontinuity of S implies the continuity of V.

The following assumption is weaker that "weakly solution-wise bounded".

$\mathcal{X}$ is said to be **weakly solution-wise bounded for the primal** if there is a non-empty bounded subset $\Omega$ of $\mathbb{R}^n$ such that for all $\xi \in \mathcal{X}$: $\Omega \cap S(\xi) \neq \phi$.

**Proposition 8:** Suppose $\mathcal{X}$ is weakly solution-wise bounded for the primal and S is upper semicontinuous on $\mathcal{X}$. Then V is continuous on $\mathcal{X}$.

**Proof:** Let $\xi = $ (p, A, b)$\in \mathcal{X}$ and let $<\xi(N)|N\in\mathbb{N}>$ be a sequence in $\mathcal{X}$ converging to $\xi$.

Sine $\mathcal{X}$ is weakly solution-wise bounded for the primal, there exists a non-empty bounded subset $\Omega$ of $\mathbb{R}^n$ such that for all $N\in \mathbb{N}$: $\Omega \cap S(\xi(N)) \neq \phi$. For each $N\in\mathbb{N}$, let $x(N)\in\Omega \cap S(\xi(N))$.

Thus, for all $N\in\mathbb{N}$, $V(\xi(N)) = p(N)^T x(N)$.

Towards a contradiction suppose that there exists $\varepsilon > 0$ such that $|V(\xi(N)) - V(\xi)| \geq \varepsilon$ infinitely often. Without loss of rigor or compromising the validity of our proof, suppose $|V(\xi(N)) - V(\xi)| \geq \varepsilon$ for all $N\in\mathbb{N}$. Such an assumption helps to economize on notation.

The sequence $<x(N)|N\in\mathbb{N}>$ in $\Omega$ is bounded since $\Omega$ is a non-empty bounded subset of $\mathbb{R}^n \times \mathbb{R}^m$. By the Bolzano-Weierstrass's theorem, $<x(N)|N\in\mathbb{N}>$ has a convergent subsequence $<x(k(N))|N\in\mathbb{N}>$ converging to some $x\in\mathbb{R}^n$. Thus, for all $N\in\mathbb{N}$, $V(\xi(k(N))) = p(k(N))^T x(k(N))$. Further, $|V(\xi(k(N))) - V(\xi)| \geq \varepsilon$ for all $N\in\mathbb{N}$.

Since S is upper semicontinuous, $<x(k(N))|N\in\mathbb{N}>$ converges to x and $x(k(N))\in S(\xi(k(N)))$ for all $N\in\mathbb{N}$ implies $x\in S(\xi)$.

Thus, $V(\xi) = p^T x = \lim_{N\to\infty} p(k(N))^T x(k(N)) = \lim_{N\to\infty} V(\xi(k(N)))$, contradicting $|V(\xi(k(N))) - V(\xi)| \geq \varepsilon$ for all $N\in\mathbb{N}$ and thereby proving the proposition. Q.E.D.

### 6. Lower-semicontinuity of solution correspondence:

Compared to upper-semicontinuity, establishing the lower-semicontinuity of the solution correspondence is far more complicated.

The following property can be found in Meyer (1979).

$\xi \in \mathbb{R}^n \times \mathbb{R}^{m\times n} \times \mathbb{R}^m$ is said to be **singleton-solvable** if $S(\xi)$ is a singleton.

The following result-in a considerably less general framework- is available in Meyer (1979).

**Proposition 9:** Suppose $\mathcal{X}$ is weakly solution-wise bounded for the primal and S is upper semicontinuous on $\mathcal{X}$. If $\xi\in\mathcal{X}$ is singleton-solvable, then S is lower- semicontinuous <u>at</u> $\xi$.

**Proof:** Since $\xi$ is singleton-solvable, there exists $x\in\mathbb{R}^n$, such that $S(\xi) = \{x\}$. Let $\xi = $ (p, A, b).

Let $<\xi(N)|N\in\mathbb{N}>$ be a sequence in $\mathcal{X}$ converging to $\xi$.

Sine $\mathcal{X}$ is weakly solution-wise bounded for the primal, there exists a non-empty bounded subset $\Omega$ of $\mathbb{R}^n$ such that for all $N\in \mathbb{N}$: $\Omega \cap S(\xi(N)) \neq \phi$. For each $N\in\mathbb{N}$, let $x(N)\in \Omega \cap S(\xi(N))$ and towards a contradiction suppose that there exists $\varepsilon > 0$ such that $\|x(N)-x\| \geq \varepsilon$ infinitely often. Without loss of rigor or compromising the validity of our proof, suppose $\|x(N) - x\| \geq \varepsilon$ for all $N\in\mathbb{N}$. Such an assumption helps to economize on notation.

The sequence $\langle x(N) | N \in \mathbb{N} \rangle$ in $\Omega$ is bounded since $\Omega$ is a non-empty bounded subset of $\mathbb{R}^n \times \mathbb{R}^m$.

By the Bolzano-Weierstrass's theorem, $\langle x(N) | N \in \mathbb{N} \rangle$ has a convergent subsequence $\langle x(k(N)) | N \in \mathbb{N} \rangle$.

Since $\|x(N) - x\| \geq \varepsilon$ for all $N \in \mathbb{N}$, it must be the case that $\|x(k(N)) - x\| \geq \varepsilon$ for all $N \in \mathbb{N}$.

Since $x(k(N)) \in S(\xi(k(N)))$ for all $N \in \mathbb{N}$ and $\langle \xi(k(N)) | N \in \mathbb{N} \rangle$ converges to $\xi$, by the upper-semicontinuity of S, $\lim_{N \to \infty} x(k(N)) \in S(\xi) = \{x\}$, thereby contradicting $\|x(k(N)) - x\| \geq \varepsilon$ for all $N \in \mathbb{N}$ and proving the proposition. Q.E.D.

**Note:** Example 1 shows how lower-semicontinuity is violated in the absence of the assumption "singleton-solvable", unless there is some other appropriate assumption to make up for its absence.

The natural question that arises is for the solution correspondence, is there any alternative to the combination of "weakly solution-wise bounded for the primal" property and upper-semicontinuity that yields lower-semicontinuity? It turns out that the requirement of the LP problem being "regular" is a sufficient alternative to the two conditions.

**Proposition 10:** Suppose $(p, A, b) \in \mathcal{X}$ is singleton-solvable and regular. Then, S is lower semicontinuous at $(p, A, b)$.

**Proof:** Suppose $(p, A, b) \in \mathcal{X}$ is singleton-solvable and regular. Then, $S(p, A, b) = S^*(p, A, b)$ is a singleton say $\{x\}$.

Thus, the array of columns $\langle A^j | x_j > 0 \rangle$ are linearly independent. Without loss of generality suppose $A = [B|E]$ where the columns of the sub-matrix B correspond to the coordinates in $\{j | x_j > 0\}$.

Thus $x = \begin{pmatrix} x_B \\ 0 \end{pmatrix} = \begin{pmatrix} B^- b \\ 0 \end{pmatrix}$.

Since $(p, A, b)$ is regular it must be the case that $(p_B)^T B^- E \gg 0$.

Let $y^T = (p_B)^T B^-$. Thus, $y^T B = (p_B)^T$ and $y^T E \gg (p_E)^T$. Thus, $y^T A \geq p^T$ and $(y^T A - p^T) x = 0$.

Let $\langle \xi(N) | N \in \mathbb{N} \rangle$ be a sequence in $\mathcal{X}$ converging to $(p, A, b)$.

Thus, there exists $N_1 \in \mathbb{N}$ such that for all $N \geq N_1$, the columns of $B(N) = \{A^j(N) | x_j > 0\}$ are linearly independent and $x(N) = \begin{pmatrix} B(N)^- b(N) \\ 0 \end{pmatrix}$ satisfies $B(N)^- b(N) \gg 0$. Further, if $y^T(N) = (p_B(N))^T B(N)^-$ then $y^T(N) B(N) = p_B(N)^T$ and $y^T(N) E(N) \gg 0$, where $E(N) = \langle A^j(N) | x_j = 0 \rangle$.

Thus, for all $N \geq N_1$, the pair $(x(N), y(N))$ satisfy the KKT conditions for $\xi(N)$ and hence $x(N) \in S(\xi(N))$.

Let $x(N) \in S(\xi(N))$ for $N < N_1$.

Clearly, $\lim_{N \to \infty} x(N) = \begin{pmatrix} (B^T B)^{-1} B^T b \\ 0 \end{pmatrix} = x$.

Thus, S is lower semicontinuous at (p, A, b). Q.E.D.

An alternative way of guaranteeing the lower-semicontinuity is to invoke the following stronger version of regularity.

(p, A, b)$\in \mathbb{R}^n \times \mathbb{R}^{m \times n} \times \mathbb{R}^m$ is said to be **strongly regular** if <u>for all</u> $x \in S^*(p, A, b)$ and $y^T = (p_B)^T B^-$ where B is the submatrix of A whose columns are $\{A^j | x_j > 0\}$ and $p_B$ is the sub-vector of p that correspond to the columns of B, it is the case that $y^T A^j - p_j > 0$ for all $j \in \{1, \ldots, n\}$ satisfying $x_j = 0$.

It is easy to see that y in the statement of the above property must satisfy $y^T A^j = p_j$ for j whenever $x_j > 0$.

The following result allows us to dispense with "singleton solvability" in order to obtain the lower-continuity of the solution correspondence at an LP problem.

**Proposition 11:** If $X(\xi)$ is bounded and $\xi$ is strongly regular then the correspondence S is lower semicontinuous at $\xi$.

**Proof:** Let $\xi = $ (p, A, b). Since $X(\xi)$ is bounded, by theorem 3 in Lahiri (2024), $S^*(\xi) = S(\xi) \cap X(\xi)$ is a non-empty finite set and all points in $S(\xi)$ can be expressed as a convex combination of points in $S^*(\xi)$.

Let $<\xi(N) | N \in \mathbb{N}>$ be a sequence in $\mathcal{X}$ converging to (p, A, b).

It is enough to show that for all $x \in S^*(\xi)$, there exists a sequence $<x(N) | N \in \mathbb{N}>$ converging to x such that $x(N) \in S(\xi(N))$ for all $N \in \mathbb{N}$.

Let $x \in S^*(\xi)$. Let B be the submatrix of A whose columns are $\{A^j | x_j > 0\}$. Without loss of generality suppose A = (B|E).

Since $\xi$ is strongly regular for $y^T = (p_B)^T B^-$, it is the case that $y^T A^j > p_j$ for all $j \in \{1, \ldots, n\}$ satisfying $x_j = 0$. Clearly, $y^T B = (p_B)^T$.

Thus, $x = \begin{pmatrix} x_B \\ x_E \end{pmatrix} = \begin{pmatrix} x_B \\ 0 \end{pmatrix}$, $p = \begin{pmatrix} p_B \\ p_E \end{pmatrix}$ and $x_B = B^- b \gg 0$.

Since, $<(p(N), A(N), b(N)) | N \in \mathbb{N}>$ converges to (p, A, b), there exists $N_1 \in \mathbb{N}$ such that for all $N \geq N_1$, $(B(N)^T B(N))^{-1}$ exists, $x_{B(N)} = B(N)^- b(N) \gg 0$ and $y(N)^T = (p_{B(N)})^T B(N)^-$ satisfies $y(N)^T A^j(N) - p_j(N) = 0$ for all j for which $x_j > 0$ and $y(N)^T A^j(N) - p_j(N) > 0$ for all j for which $x_j = 0$.

Let $x(N) = \begin{pmatrix} x_{B(N)} \\ x_{E(N)} \end{pmatrix} = \begin{pmatrix} x_{B(N)} \\ 0 \end{pmatrix}$ for all $N \geq N_1$. Thus $\{j | x_j > 0\} = \{j | x_j(N) > 0\}$ for all $N \geq N_1$. Further $(p_{B(N)})^T x_{B(N)} = y(N)^T b(N)$ for all $N \geq N_1$. Thus, $x(N) \in S(\xi(N))$ for all $N \geq N_1$.

Let $x(N) \in S(\xi(N))$ for all $N < N_1$.

Clearly, $\lim_{N \to \infty} x(N) = \begin{pmatrix} (B^T B)^{-1} B^T b \\ 0 \end{pmatrix} = x$.

Thus, S is lower semicontinuous at $\xi$. Q.E.D.

**7. Domain with variable objective function:**

$\mathcal{X}$ is said to be **a domain with variable objective function** if there exists $A \in \mathbb{R}^{m \times n}$ and $b \in \mathbb{R}^m$, such that for all $\xi \in \mathcal{X}$, $A(\xi) = A$ and $b(\xi) = b$.

Let $(p, A, b) \in \mathcal{X}$ and without loss of generality suppose $A = (B|E)$ satisfies the following properties:

(i) the columns of A comprising the submatrix B are linearly independent;

(ii) For $x = \begin{pmatrix} x_B \\ x_E \end{pmatrix} = \begin{pmatrix} B^- b \\ 0 \end{pmatrix}$ with $x_B = B^- b$, and $p = \begin{pmatrix} p_B \\ p_E \end{pmatrix}$ it is the case that $(p_B)^T B^- E \gg p_E$.

Thus, $x \in S^*(p, A, b)$ and $y = (p_B)^T B^-$ solves $D(p, A, b)$.

Let $\Delta p$ be a non-zero n-dimensional column vector and consider $(p + \theta \Delta p, A, b)$.

So long as $(p_B + \theta \Delta p_B)^T B^- E \geq p_E + \theta \Delta p_E$, $x \in S^*(p + \theta \Delta p, A, b)$.

Thus, so long as $(p_B + \theta \Delta p_B)^T B^- E^j \geq p_j + \theta \Delta p_j$ for all $E^j$, it must be the case that $x \in S^*(p + \theta \Delta p, A, b)$.

Thus, $\theta[(\Delta p_B)^T B^- E^j - \Delta p_j] \geq p_j - (p_B)^T B^- E^j$ implies $x \in S^*(p + \theta \Delta p, A, b)$.

Thus, if $\theta \geq \max\{\frac{p_j - (p_B)^T B^- E^j}{[(\Delta p_B)^T B^- E^j - \Delta p_j]} | (\Delta p_B)^T B^- E^j - \Delta p_j > 0\}$ and $\theta \leq \min\{\frac{p_j - (p_B)^T B^- E^j}{[(\Delta p_B)^T B^- E^j - \Delta p_j]} | (\Delta p_B)^T B^- E^j - \Delta p_j < 0\}$, then $x \in S^*(p + \theta \Delta p, A, b)$

Since, $p_j - (p_B)^T B^- E^j < 0$ for all $E^j$, $\max\{\frac{p_j - (p_B)^T B^- E^j}{[(\Delta p_B)^T B^- E^j - \Delta p_j]} | (\Delta p_B)^T B^- E^j - \Delta p_j > 0\} < 0 < \min\{\frac{p_j - (p_B)^T B^- E^j}{[(\Delta p_B)^T B^- E^j - \Delta p_j]} | (\Delta p_B)^T B^- E^j - \Delta p_j < 0\}$.

Hence for all $\theta \in [\max\{\frac{p_j - (p_B)^T B^- E^j}{[(\Delta p_B)^T B^- E^j - \Delta p_j]} | (\Delta p_B)^T B^- E^j - \Delta p_j > 0\}, \min\{\frac{p_j - (p_B)^T B^- E^j}{[(\Delta p_B)^T B^- E^j - \Delta p_j]} | (\Delta p_B)^T B^- E^j - \Delta p_j < 0\}]$, $V(p + \theta \Delta p, A, b) = (p_B + \theta \Delta p_B)^T x_B = V(p, A, b) + \theta(\Delta p_B)^T x_B = V(p, A, b) + \theta(\Delta p_B)^T B^- b$.

As in the framework of in Bohm (1975), in the case of domain with variable objective function, this local linearity property of the optimal value function does not depend on whether for $\theta \neq 0$, $(p + \theta \Delta p, A, b) \in \mathcal{X}$ or not.